\documentclass{amsart}

\usepackage{cite}
\usepackage{textcomp}

\usepackage{amsmath, amssymb, amsfonts}
\usepackage{accents}
\usepackage{algorithm}
\usepackage{algorithmic}
\usepackage{color}
\usepackage{graphicx}
\usepackage{enumerate}
\usepackage{booktabs}
\usepackage{subcaption}
\usepackage{enumerate}
\usepackage{mathrsfs}
\usepackage{mathtools}
\usepackage{dsfont}

\DeclareMathOperator{\tr}{tr}
\DeclareMathOperator{\minimize}{minimize}

\DeclareMathOperator{\voi}{VoI}
\DeclareMathOperator{\E}{\mathsf{E}}
\DeclareMathOperator{\Cov}{\mathsf{cov}}

\DeclareMathOperator{\ProbM}{\mathsf{P}}
\DeclareMathOperator{\Prob}{\mathsf{p}}
\DeclareMathOperator{\ProbG}{\mathsf{s}}
\DeclareMathOperator{\ProbQ}{\mathsf{q}}

\newtheorem{definition}{Definition}
\newtheorem{lemma}{Lemma}
\newtheorem{theorem}{Theorem}

\newtheorem{remark}{Remark}

\begin{document}
\title{Value of Information in Feedback Control: Global~Optimality}

 \author{
Touraj Soleymani, John S. Baras, Sandra Hirche, and Karl H. Johansson}
\thanks{Corresponding Author: Touraj Soleymani (touraj@kth.se). Journal: \emph{IEEE Transactions on Automatic Control}.}

\maketitle
\begin{abstract}
The rate-regulation tradeoff, defined between two objective functions, one penalizing the packet rate and one the regulation cost, can express the fundamental performance bound of networked control systems. However, the characterization of the set of globally optimal solutions in this tradeoff for multi-dimensional Gauss--Markov processes has been an open problem. In the present article, we characterize a policy profile that belongs to this set without imposing any restrictions on the information structure or the policy structure. We prove that such a policy profile consists of a symmetric threshold triggering policy based on the value of information and a certainty-equivalent control policy based on a non-Gaussian linear estimator. These policies are deterministic and can be designed separately. Besides, we provide a global optimality analysis for the value of information $\voi_k$, a semantic metric that emerges from the rate-regulation tradeoff as the difference between the benefit and the cost of a data packet. We prove that it is globally optimal that a data packet containing sensory information at time $k$ be transmitted to the controller only if $\voi_k$ becomes nonnegative. These results have important implications in the areas of communication and control.

\smallskip
\noindent \textbf{Keywords.}
decision policies, globally optimal solutions, networked control systems, rate-regulation tradeoff, semantic communications, semantic metrics, value of information.
\end{abstract}

\section{Introduction}
The rate-regulation tradeoff, defined between two objective functions, one penalizing the packet rate and one the regulation cost, can express the \emph{fundamental performance bound} of networked control systems. Such a tradeoff naturally leads to the adoption of an event trigger that is collocated with the sensor and of a controller that is collocated with the actuator as the distributed decision makers, and is formulated as a stochastic optimization problem over the space of causal decision policy profiles. Unfortunately, this optimization problem for the joint design of the event trigger and the controller is in general intractable~\cite{wu2013, ramesh2013}. Despite lack of a general theory for coping with this difficulty, our goal here is to find a globally optimal solution in the rate-regulation tradeoff, and provide a global optimality analysis for the value of information, a quantity that emerges from the rate-regulation tradeoff and systematically captures the semantics of data packets by taking into account their potential impacts. We previously argued in~\cite{voi} that the value of information as a semantic metric determines the \emph{right piece of information}, a concept that is not defined in classical data communication, while it is crucial to the development of future communication networks. In this respect, the goal we pursue here not only is interesting on its own from a theoretical perspective, but, if achieved, has important implications in the areas of communication and control. 

In what follows, we first review and categorize the previous studies on networked systems that are closely related to our work, and then provide an overview of our~results.

\subsection{Related Work}
There exist a number of studies that have explored a tradeoff between the packet rate and the mean-square error, and characterized the optimal triggering policy~\cite{imer2010, lipsa2011, molin2017, chakravorty2016, rabi2012, guo2021-IT}. The intrinsic difficulty in these studies is due to a non-classical information structure, which complicates the derivation of the optimal triggering policy. Notably, Imer and Ba{\c{s}}ar \cite{imer2010} studied the optimal event-triggered estimation of a scalar Gauss--Markov process based on dynamic programming by assuming that the triggering policy is symmetric threshold, and derived the optimal threshold value of the policy. Lipsa and Martins~\cite{lipsa2011} analyzed the optimal event-triggered estimation of a scalar Gauss--Markov process based on majorization theory, and proved that the optimal triggering policy is symmetric threshold. Molin and Hirche~\cite{molin2017} studied the convergence properties of an iterative algorithm for the optimal event-triggered estimation of a scalar Markov process with symmetric noise distribution, and found a result coinciding with that in~\cite{lipsa2011}. Chakravorty and Mahajan~\cite{chakravorty2016} addressed the optimal event-triggered estimation of a scalar autoregressive Markov process with symmetric noise distribution based on renewal theory, and proved that the optimal triggering policy remains symmetric threshold. In addition, Rabi~\emph{et~al.}~\cite{rabi2012} formulated the optimal event-triggered estimation of the scalar Ornstein--Uhlenbeck process as an optimal multiple stopping time problem by assuming that the estimator is linear, and showed that the optimal triggering policy is symmetric threshold. Guo and Kostina \cite{guo2021-IT} also contributed to this area by studying the optimal event-triggered estimation of the scalar Ornstein--Uhlenbeck process without any assumption on the estimator, and obtained a similar result as~in~\cite{rabi2012}.

Aside from the above line of research, several works have investigated optimal event-triggered estimation when the triggering policy is fixed~\cite{sijs2012, wu2013, he2018, han2015}. The main challenge in these works is to find a procedure for dealing with a signaling effect, which can cause a nonlinearity in the structure of the optimal estimator. To that end, Sijs and Lazar \cite{sijs2012} used a sum of Gaussian approximation, and developed an estimator that has an asymptotically bounded estimation error covariance for a Gauss--Markov process subject to a fixed deterministic triggering policy. Wu~\emph{et~al.}~\cite{wu2013} used a Gaussian approximation, and found a suboptimal estimator for a Gauss--Markov process subject to a fixed deterministic threshold triggering policy. He~\emph{et~al.}~\cite{he2018} took one step further, and adopted the generalized closed skew normal distribution to characterize the optimal estimator for a Gauss--Markov process subject to a similar triggering policy. Han~\emph{et~al.}~\cite{han2015} also took advantage of a fixed stochastic triggering policy that preserves the Gaussianity of the conditional distribution, and obtained the optimal estimator for a Gauss--Markov process.

Furthermore, several works have investigated optimal event-triggered control when the triggering policy is fixed~\cite{ramesh2013, molin2013, demirel2018}. Note that this problem is more complicated than the estimation counterpart because of a dual effect, which can lead to a coupling between estimation and control. In this context, Molin and Hirche~\cite{molin2013} studied the optimal event-triggered control of a Gauss--Markov process, and showed that the optimal control policy is certainty equivalent when the triggering policy is reparametrizable in terms of primitive random variables. Ramesh~\emph{et~al.}~\cite{ramesh2013} studied the dual effect in the optimal event-triggered control of a Gauss--Markov process, and proved that the dual effect in general exists. They also proved that the certainty equivalence principle holds if and only if the triggering policy is independent of the control policy. Later, Demirel~\emph{et~al.}~\cite{demirel2018} addressed the optimal event-triggered control of a Gauss--Markov process by adopting a stochastic triggering policy that preserves the Gaussianity of the conditional distribution, and showed that the optimal control policy remains certainty equivalent.

On the contrary to the above vein of research, there exist a few studies that have considered a tradeoff between the packet rate and the trace of variance~\cite{leong2017, leong2018}. In this case, one instead of an observation-based triggering policy, i.e., the type used in \cite{voi, imer2010, lipsa2011, molin2017, chakravorty2016, rabi2012, guo2021-IT, sijs2012, wu2013, he2018, han2015, ramesh2013, molin2013, demirel2018}, searches for a variance-based triggering policy. These studies are somehow related to sensor scheduling, which dates back to a few decades ago~\cite{kushner}. Previously, Kushner~\cite{kushner} studied the optimal control of a Gauss--Markov process subject to a limited number of observations, and found the optimal triggering policy that does not depend on the observations. Recently, Leong~\emph{et~al.}~\cite{leong2017, leong2018} addressed the optimal variance-based event-triggered estimation of a Gauss--Markov process, and showed that the optimal triggering policy is a threshold policy that can be expressed in terms of the estimation error covariance. Note that when a variance-based triggering policy is used, the certainty equivalence principle simply holds~\cite{leong2017-ifac}. Nevertheless, variance-based triggering policies are generally outperformed by observation-based triggering policies, as they do not take advantage of realized sensory information.

Moreover, there exist a few studies that have considered a tradeoff between the bit rate and the mean-square error in a causal setting\cite{witsenhausen1979, walrand1983, borkar2001}. In this case, one instead of a triggering policy searches for a quantization policy. In particular, Witsenhausen~\cite{witsenhausen1979} addressed the sequential coding of a discrete-time $k$-th order Markov process over a finite time horizon, and showed that the optimal code depends on the last $k$ process states and the current decoder state. Walrand and Varaiya~\cite{walrand1983} investigated the sequential coding of a discrete-time finite-state Markov process over a noisy channel with feedback, and showed that there exists a separation in the design of the encoder and the decoder through the conditional distribution. Borkar~\emph{et al.}~\cite{borkar2001} also studied the sequential coding of a discrete-time Markov process without fixing the quantization levels, and provided a procedure based on dynamic programming for the computation of the optimal partition. Later, Y\"{u}ksel \cite{yuksel2013} extended the above results to optimal control, and showed that for a Gauss--Markov process the globally optimal quantization policy is predictive and the globally optimal control policy is certainty equivalent. Note that all these studies assume that quantized sensory information is transmitted in a~periodic~way.

\subsection{Overview and Outline}
Despite a considerable body of research in the area of networked systems, the characterization of the set of globally optimal solutions in the rate-regulation tradeoff, as described above, for multi-dimensional Gauss--Markov processes has been an open problem. In the present article, we characterize for the first time a policy profile that belongs to this set without imposing any restrictions on the information structure or the policy structure. We prove that such a policy profile consists of a symmetric threshold triggering policy and a certainty-equivalent control policy. More specifically, we show that the rate-regulation tradeoff attains a globally optimal solution of the form $(\pi^{\star},\mu^{\star}) = ( \{\mathds{1}_{\voi_k \geq 0}\}_{k=0}^{N}, \{ - L_k \hat{x}_k \big\}_{k=0}^{N} )$, where $\mathds{1}_{\voi_k \geq 0}$ denotes the indicator function of ${\voi_k \geq 0}$, $\voi_k$ is the value of information, $L_k$ is the linear-quadratic-regulator gain, and $\hat{x}_k$ is the minimum mean-square-error state estimate at the controller. Clearly, our study is different from the studies in~\cite{imer2010, lipsa2011, molin2017, chakravorty2016, rabi2012, guo2021-IT}, where the results apply to the estimation of scalar processes. Here, the results apply to the control of multi-dimensional Gauss--Markov processes. Our study is also different from the studies in~\cite{sijs2012, wu2013, han2015, he2018, molin2013, ramesh2013, demirel2018}, where an estimation policy or a control policy is derived when the triggering policy is fixed and subject to some conditions. Here, we search for a globally optimal triggering policy and a globally optimal control policy jointly and without any restrictions. Finally, our study differs from the studies in~\cite{leong2017, leong2018, leong2017-ifac, witsenhausen1979, walrand1983, borkar2001, yuksel2013}, where a variance-based triggering policy or a quantization policy is derived. Here, we are particularly interested in observation-based triggering~policies.

Besides, in this article, we provide for the first time a global optimality analysis for the value of information $\voi_k$, which in fact measures the difference between the benefit and the cost of a data packet. We prove that it is globally optimal that a data packet containing sensory information at time $k$ be transmitted to the controller only if $\voi_k$ becomes nonnegative. Using backward induction in~\cite{voi}, we quantified and approximated the value of information for multi-dimensional Gauss--Markov processes at a Nash equilibrium, where neither decision maker has a unilateral incentive to change its policy. However, a question that was not addressed there is whether this equilibrium is globally optimal. The importance of this question cannot be overstated, as the rate-regulation tradeoff might admit other Nash equilibria with better performance. We address this question in the present article by developing new techniques, and prove that the previously characterized Nash equilibrium has zero optimality gap. Throughout our analysis, we will use the existence result and some of the mathematical derivations of~\cite{voi}.

The article is organized in the following way. We formulate the rate-regulation tradeoff in Section~\ref{sec2}, and present our main result in Section~\ref{sec3}. Finally, we conclude the article in Section~\ref{sec4}.

\subsection{Preliminaries}
In the sequel, the sets of real numbers and non-negative integers are denoted by $\mathbb{R}$ and $\mathbb{N}$, respectively. For $x,y \in \mathbb{N}$ and $x \leq y$, the set $\mathbb{N}_{[x,y]}$ denotes $\{z \in \mathbb{N} | x \leq z \leq y\}$. The sequence of vectors $x_{0}, \dots, x_{k}$ is represented by $\mathbf{x}_k$. For matrices $X$ and $Y$, the relations $X \succ 0$ and $Y \succeq 0$ denote that $X$ and $Y$ are positive definite and positive semi-definite, respectively. The indicator function of a subset $\mathcal{A}$ of a set $\mathcal{X}$ is denoted by $\mathds{1}_\mathcal{A}:\mathcal{X} \to \{0,1\}$. The symmetric decreasing rearrangement of a Borel measurable function $f(x)$ vanishing at infinity is represented by $f^*(x)$. The probability measure of a random variable $x$ is represented by $\mathsf{P}(x)$, its probability density or probability mass function by $\Prob(x)$, and its expected value and covariance by $\E[x]$ and $\Cov[x]$, respectively.

\begin{definition}(Stochastic kernels)
Let $(\mathcal{X},\mathcal{B}_{\mathcal{X}})$ and $(\mathcal{Y},\mathcal{B}_{\mathcal{Y}})$ be two measurable spaces. A Borel measurable stochastic kernel $\ProbM: \mathcal{B}_{\mathcal{Y}} \times  \mathcal{X} \to [0,1]$ is a mapping such that $\mathcal{A} \mapsto \ProbM( \mathcal{A} | x)$ is a probability measure on $(\mathcal{Y},\mathcal{B}_{\mathcal{Y}})$ for any $x \in \mathcal{X}$, and $x \mapsto \ProbM(\mathcal{A}| x )$ is a Borel measurable function for any $\mathcal{A} \in \mathcal{B}_{\mathcal{Y}}$.
\end{definition}

\begin{definition}(Globally optimal solutions)
For a given team game with two decision makers, let $\gamma^1 \in \mathcal{G}^1$ and $\gamma^2 \in \mathcal{G}^2$ be the decision policies of the decision makers, where $\mathcal{G}^1$ and $\mathcal{G}^2$ are the sets of admissible policies, and $L(\gamma^1,\gamma^2)$ be the associated loss function. A policy profile $(\gamma^{1\star},\gamma^{2\star})$ is globally optimal if
\begin{align*}
	L(\gamma^{1\star},\gamma^{2\star}) \leq L(\gamma^{1},\gamma^{2}), \ \text{for all } \gamma^1 \in \mathcal{G}^1, \gamma^2 \in \mathcal{G}^2.
\end{align*}
Note that globally optimal solutions express a stronger solution concept than Nash equilibria.
\end{definition}

\section{Rate-Regulation Tradeoff}\label{sec2}
Consider a networked control system in its basic form. The dynamics of the underlying process is given by the discrete-time state and output equations
\begin{align}
	x_{k+1} &= A_k x_k + B_k u_k + w_k,\label{c1:eq:sys}\\[1.75\jot]
	y_k &= C_k x_k + v_k,\label{c1:eq:output-i}
\end{align}
for $k \in \mathbb{N}_{[0,N]}$ with initial condition $x_0$, where $x_k \in \mathbb{R}^n$ is the state of the process, $A_k \in \mathbb{R}^{n \times n}$ is the state matrix, $B_k \in \mathbb{R}^{n \times m}$ is the input matrix, $u_k \in \mathbb{R}^m$ is the control input applied by an actuator and decided by a controller that is collocated with the actuator, $w_k \in \mathbb{R}^n$ is a Gaussian white noise with zero mean and covariance $W_k \succ 0$, $y_k \in \mathbb{R}^{p}$ is the output of the process observed by a sensor, $C_k \in \mathbb{R}^{p \times n}$ is the output matrix, $v_k \in \mathbb{R}^{p}$ is a Gaussian white noise with zero mean and covariance $V_k \succ 0$, and $N \in \mathbb{N}$ is a finite time horizon. It is assumed that $x_0$ is a Gaussian vector with mean $m_0$ and covariance $M_0$,  and that $x_0$, $w_k$, and $v_k$ are mutually independent for all $k \in \mathbb{N}_{[0,N]}$. The feedback control loop is closed via a reliable but costly communication channel, and the sensory information in this channel is carried in the form of data packets subject to one-step delay. Let $a_k$ and $b_k$ represent the input and the output of the channel at time $k$, respectively. Then, we have
\begin{align}\label{c1:eq:etm1}
b_{k+1} = \left\{
  \begin{array}{l l}
     a_k, & \ \text{if} \ \delta_k =1, \\
     \varnothing, & \ \text{otherwise},
  \end{array} \right.
\end{align}
for $k \in \mathbb{N}_{[0,N]}$ with $b_0 = \varnothing$, where $\delta_k \in \{0,1\}$ is the transmission decision decided by an event trigger that is collocated with the sensor. It is assumed that the data packet that can be transmitted at time $k$ contains the minimum mean-square-error state estimate at the event trigger at time $k$, and that the quantization error is negligible.

The event trigger and the controller, as two distributed decision makers, make their decisions based on their causal information sets, which are given by $\mathcal{I}^e_k := \{ y_t, b_t, \delta_{s}, u_{s} | t \in \mathbb{N}_{[0,k]}, s \in \mathbb{N}_{[0,k-1]} \}$ and $	\mathcal{I}^c_k := \{ b_t, \delta_{s}, u_{s} | t \in \mathbb{N}_{[0,k]}, s \in \mathbb{N}_{[0,k-1]} \}$, respectively. We say that a triggering policy $\pi$ and a control policy $\mu$ are admissible if $\pi = \{\ProbM(\delta_k | \mathcal{I}^e_k) \}_{k=0}^{N}$ and $\mu = \{\ProbM(u_k | \mathcal{I}^c_k) \}_{k=0}^{N}$, where $\ProbM(\delta_k | \mathcal{I}^e_k)$ and $\ProbM(u_k | \mathcal{I}^c_k)$ are Borel measurable stochastic kernels. We represent the sets of admissible triggering policies and admissible control policies by $\mathcal{P}$ and $\mathcal{M}$, respectively.

Our goal in this study is to find a globally optimal solution $(\pi^\star,\mu^\star)$ to the following stochastic optimization problem:
\begin{align}\label{eq:main_problem1}
	\underset{\pi \in \mathcal{P},\mu \in \mathcal{M}}{\minimize} \  \Phi(\pi,\mu) := (1-\lambda) R(\pi,\mu) + \lambda J(\pi,\mu),
\end{align}
for the tradeoff multiplier $\lambda \in (0,1)$ and
\begin{align}
R(\pi,\mu) &:= \textstyle \frac{1}{N+1} \E\Big[\sum_{k=0}^{N} \ell_k \delta_k \Big],\label{eq:rate-measure}\\[1\jot]
J(\pi,\mu) &:= \textstyle \frac{1}{N+1} \E \Big [\textstyle\sum_{k=0}^{N} x_{k+1}^T Q_{k+1} x_{k+1} + u_k^T R_{k} u_k \Big],\label{eq:control-measure}
\end{align}
where $\ell_k \geq 0$ is a weighting coefficient and $Q_k \succeq 0$ and $R_k \succ 0$ are weighting matrices.

\begin{remark}
The optimization problem in (\ref{eq:main_problem1}) formulates the rate-regulation tradeoff between the packet rate and the regulation cost for multi-dimensional Gauss--Markov processes. Note that the set of globally optimal solutions in this tradeoff cannot be empty following our results in \cite{voi}, where the existence of a Nash equilibrium is proved. In the sequel, we in fact investigate the optimality gap of this very equilibrium. Our study focuses on the soft-constraint version of the rate-regulation tradeoff, where the packet rate appears in the loss function. The hard-constraint version of the rate-regulation tradeoff, where the packet rate appears as a constraint, attains the same solutions as long as there exists an associated Lagrange multiplier.
\end{remark}

\section{Global Optimality Analysis of\\the Value of Information}\label{sec3}
The main result of this article is provided in this section. We first introduce two distinct value functions from the perspectives of the event trigger and the controller, and then provide the general formula of the value of information.

\begin{definition}[Value functions]
The value functions $V^e_k(\mathcal{I}^e_k)$ and $V^c_k(\mathcal{I}^c_k)$ are defined as
\begin{align}
	V^e_k(\mathcal{I}^e_k) :=& \min_{\pi \in \mathcal{P} : \mu = \mu^\star}\E \Big[ \textstyle \sum_{t=k}^{N} \theta_t \delta_t + \varsigma_{t+1} \Big| \mathcal{I}^e_k \Big],\label{eq:Ve-def}\\[1\jot]
	V^c_k(\mathcal{I}^c_k) :=& \min_{\mu \in \mathcal{M}: \pi = \pi^\star}\E \Big[ \textstyle \sum_{t=k}^{N} \theta_{t-1} \delta_{t-1} + \varsigma_{t} \Big| \mathcal{I}^c_k \Big],\label{eq:Vc-def}
\end{align}
for $k \in \mathbb{N}_{[0,N]}$ given a policy profile $(\pi^\star,\mu^\star)$, where $\theta_k = \ell_k (1-\lambda)/\lambda$ and $\varsigma_k =  (u_k + (B_k^T S_{k+1} B_k + R_k)^{-1} B_k^T S_{k+1} A_k x_k )^T (B_k^T S_{k+1} B_k + R_k) (u_k + (B_k^T S_{k+1} B_k + R_k)^{-1} B_k^T S_{k+1} A_k x_k )$ with the exception of $\theta_{-1} = 0$ and $\varsigma_{N+1} = 0$, and $S_k \succeq 0$ obeys the algebraic Riccati equation
\begin{equation}\label{eq:riccati}
\begin{aligned}
S_k &= Q_k + A_k^T S_{k+1} A_k - A_k^T S_{k+1} B_k \\[1.75\jot]
	&\qquad \times (B_k^T S_{k+1} B_k + R_k)^{-1} B_k^T S_{k+1} A_k,
\end{aligned}
\end{equation}
for $k \in \mathbb{N}_{[0,N]}$ with initial condition $S_{N+1} = Q_{N+1}$ and with the exception of $S_{N+2} = 0$.
\end{definition}

\begin{definition}[Value of Information]\label{def:voi}
The value of information at time $k$ is defined as the variation in the value function $V^e_k(\mathcal{I}^e_k)$ with respect to the sensory information $a_k$ that can be communicated to the controller at time~$k$, i.e.,
\begin{align}\label{eq:voi-def}
\voi_k := V^e_k(\mathcal{I}^e_k)|_{\delta_k = 0} - V^e_k(\mathcal{I}^e_k)|_{\delta_k = 1},
\end{align}
where $V^e_k(\mathcal{I}^e_k)|_{\delta_k}$ denotes the value function $V^e_k(\mathcal{I}^e_k)$ when the transmission decision $\delta_k$ is enforced.
\end{definition}

Let $\check{x}_k := \E[x_k | \mathcal{I}^e_k]$ and $\hat{x}_k := \E[x_k | \mathcal{I}^c_k]$ denote the minimum mean-square-error state estimates at the event trigger and the controller, respectively. In addition, let us define the estimation error from the perspective of the event trigger $\check{e}_k := x_k - \E[x_k | \mathcal{I}^e_k]$, the estimation error from the perspective of the controller $\hat{e}_k := x_k - \E[x_k | \mathcal{I}^c_k]$, and the estimation mismatch $\tilde{e}_k := \E[x_k | \mathcal{I}^e_k] - \E[x_k | \mathcal{I}^c_k]$. The next theorem states our main result on the characterization of a globally optimal solution in the rate-regulation tradeoff.

\begin{theorem}\label{thm:1}
The rate-regulation tradeoff attains a globally optimal solution $(\pi^{\star},\mu^{\star})$ such that
\begin{align}\label{eq:opt-profile}
(\pi^{\star},\mu^{\star}) = \Big( \big\{\mathds{1}_{\voi_k \geq 0} \big\}_{k=0}^{N}, \big\{ - L_k \hat{x}_k \big\}_{k=0}^{N} \Big),
\end{align}
with
\begin{align}
	\voi_k &= \tilde{e}_k^T A_k^T \Gamma_{k+1} A_k \tilde{e}_k- \theta_k + \varrho_k,\label{eq:voi-stale}\\[1.75\jot]
	\hat{x}_{k+1} &= A_k \hat{x}_k + B_k u_k + \delta_k A_k \tilde{e}_k,\label{eq:linear-filter}
\end{align}
for $k \in \mathbb{N}_{[0,N]}$, where $L_k = (B_k^T S_{k+1} B_k + R_k)^{-1} B_k^T S_{k+1} A_k$ is the control gain, $\Gamma_k = A_k^T S_{k+1} B_k (B_k^T S_{k+1} B_k + R_k)^{-1} B_k^T S_{k+1} A_k$ is a weighting matrix, $\varrho_k = \E[V^e_{k+1}(\mathcal{I}^e_{k+1})|$ $\mathcal{I}^e_k, \delta_k = 0] - \E[V^e_{k+1}(\mathcal{I}^e_{k+1})|\mathcal{I}^e_k, \delta_k = 1]$ is a symmetric function of $\tilde{e}_k$, and $\hat{x}_0 = m_0$ is the initial condition.
\end{theorem}

\begin{remark}
The globally optimal solution $(\pi^{\star},\mu^{\star})$ in (\ref{eq:opt-profile}) consists of a symmetric threshold triggering policy based on the value of information and a certainty-equivalent control policy based on a non-Gaussian linear state estimator. This result is important as it shows that the characterized Nash equilibrium in~\cite{voi} has zero optimality gap. Observe that the decision policies $\pi^{\star}$ and $\mu^{\star}$ are deterministic, implying that randomization does not improve the system performance, and that they can be designed separately. Moreover, note that $\voi_k(\mathcal{I}^e_k)$ in (\ref{eq:voi-stale}), which is a symmetric function of the estimation mismatch $\tilde{e}_k$, measures the difference between the benefit of transmitting a data packet, i.e., $\tilde{e}_k^T A_k^T \Gamma_{k+1} A_k \tilde{e}_k + \varrho_k$, and its associated cost, i.e., $\theta_k$. This means that it is globally optimal that a data packet containing the sensory information $\check{x}_k$ be transmitted to the controller only if its benefit surpasses its cost, i.e., $\voi_k \geq 0$. Furthermore, note that the state estimate $\hat{x}_k$ in (\ref{eq:linear-filter}) obeys a linear recursive equation with no residual $\imath_k := A_k \E[\hat{e}_k | \mathcal{I}^c_k,\delta_k=0]$ (see Lemma~\ref{lemma:estimator-at-controller} in the Appendix for the general equation of the optimal estimator at the controller). This implies that the controller's inference about the state of the process when no data packet is delivered has no contribution from the minimum mean-square-error perspective. Finally, we remark that at the globally optimal solution $(\pi^\star,\mu^\star)$ the transmission of the state estimate $\check{x}_k$ is equivalent to that of the estimation mismatch $\tilde{e}_k$, whose magnitude is comparatively smaller.
\end{remark}

\begin{proof}
Let $(\pi^o,\mu^o)$ denote a policy profile in the set of globally optimal solutions. As we said earlier, this set cannot be empty. We prove that the policy profile $(\pi^\star,\mu^\star)$ in the claim is globally optimal by showing that $\Phi(\pi^\star,\mu^\star)$ cannot be greater than $\Phi(\pi^o,\mu^o)$. Our proof is structured in the following way. We first find an innovation-based triggering policy $\sigma$ such that $\Phi(\sigma,\mu^o) = \Phi(\pi^o,\mu^o)$. Then, we derive a certainty-equivalent control policy $\xi$ such that $\Phi(\sigma,\xi) \leq \Phi(\sigma,\mu^o)$. Afterwards, we construct a symmetric triggering policy $\omega$ such that $\Phi(\omega,\xi) \leq \Phi(\sigma,\xi)$. Finally, we show that for the policy profile in the claim we have $\Phi(\pi^\star,\mu^\star) \leq \Phi(\omega,\xi)$. Throughout our analysis, without loss of generality, we assume that $m_0 = 0$. Similar arguments can be made for $m_0 \neq 0$ following a coordinate transformation.

In the first step, we will show that, given the control policy $\mu^o$, we can find an innovation-based triggering policy $\sigma$ that is equivalent to the triggering policy $\pi^o$. Note that the innovation $\nu_k := y_k - C_k \E[ x_k | \mathcal{I}^e_{k-1}]$ is a white Gaussian noise with zero mean and covariance $N_k = C_k M_k C_k^T + V_k$, where $M_k = \Cov[x_k | \mathcal{I}^e_{k-1}]$. From this definition, we have $\mathbf{y}_k = \boldsymbol{\nu}_k + E_k \check{\mathbf{x}}_{k-1} + F_k \mathbf{u}_{k-1}$, where $E_k$ and $F_k$ are matrices of proper dimensions. By Lemma~\ref{lemma:Kalmanfilter}, we have $\check{\mathbf{x}}_k = G_k \boldsymbol{\nu}_k + H_k \mathbf{u}_{k-1}$, where $G_k$ and $H_k$ are matrices of proper dimensions. In addition, from (\ref{c1:eq:etm1}), we know that $\mathbf{b}_k$ is a function of $\check{\mathbf{x}}_{k-1}$ and $\boldsymbol{\delta}_{k-1}$. As a result, it is possible to write
\begin{align*}
\Prob_{\pi^o}(\delta_k | \mathcal{I}^e_k) &= \Prob_{\pi^o}(\delta_k | \boldsymbol{\nu}_k, \boldsymbol{\delta}_{k-1}, \mathbf{u}_{k-1}),\\[2.25\jot]
\Prob_{\mu^o}(u_k | \mathcal{I}^c_k) &= \Prob_{\mu^o}(u_k | \boldsymbol{\nu}_{k-1}, \boldsymbol{\delta}_{k-1}, \mathbf{u}_{k-1}).
\end{align*}
Accordingly, any realizations of $\delta_k$ and $u_k$ can be expressed as $\delta_k = \delta_k \big( \eta_k;\boldsymbol{\nu}_k, \boldsymbol{\delta}_{k-1},$ $\mathbf{u}_{k-1} \big)$ and  $u_k = u_k \big( \zeta_k; \boldsymbol{\nu}_{k-1}, \boldsymbol{\delta}_{k-1}, \mathbf{u}_{k-1} \big)$, respectively, where $\eta_k$ and $\zeta_k$ represent random variables that are independent of any other variables. Hence, it is possible to recursively construct $\sigma$ with $\Prob_{\sigma}(\delta_k | \boldsymbol{\nu}_k, \boldsymbol{\delta}_{k-1}, \boldsymbol{\zeta}_{k-1})$ such that it is equivalent to $\Prob_{\pi^o}(\delta_k | \mathcal{I}^e_k)$. This proves that $\Phi(\sigma,\mu^o) = \Phi(\pi^o,\mu^o)$. Note that although the triggering policy $\sigma$ has been constructed associated with the control policy $\mu^o$, it now depends only on $\boldsymbol{\nu}_k$, $\boldsymbol{\delta}_{k-1}$, and $\boldsymbol{\zeta}_{k-1}$ at each time $k$.

In the second step, given the triggering policy $\sigma$, we will search for an optimal control policy $\xi$, and prove that $\xi$ is certainty equivalent. Using (\ref{c1:eq:sys}) and (\ref{eq:riccati}), we can derive the following identities:
\begin{align}
\begin{split}\label{eq:identity1-1}
	&x_{k+1}^T S_{k+1} x_{k+1} = (A_k x_k + B_k u_k + w_k)^T\\[1.95\jot]
	&\quad\qquad \qquad \qquad \quad \times S_{k+1} 	(A_k x_k + B_k u_k + w_k),
\end{split}\\[1.75\jot]
\begin{split}\label{eq:identity1-2}
	&x_k^T S_k x_k = x_k^T \big(Q_k + A_k^T S_{k+1} A_k\\[1.95\jot]
	&\quad\qquad \qquad \qquad - L_k ^T (B_k^T S_{k+1} B_k + R_k) L_k\big) x_k,
\end{split}\\[1.75\jot]
\begin{split}\label{eq:identity1-3}
	&x_{N+1}^T S_{N+1} x_{N+1} - x_0^T S_0 x_0\\[1\jot]
	&\quad \qquad =  \textstyle \sum_{k=0}^{N} x_{k+1}^T S_{k+1} x_{k+1} - \sum_{k=0}^{N} x_k^T S_k x_k.
\end{split}
\end{align}
Then, incorporating the identities (\ref{eq:identity1-1}) and (\ref{eq:identity1-2}) into the identity (\ref{eq:identity1-3}), taking the expectation of both sides of (\ref{eq:identity1-3}), and using the facts that $w_k$ is independent of $x_k$ and $u_k$ and that the terms $x_0^T S_0 x_0$ and $w_k^T S_{k+1} w_k$ are independent of the decision policies, we find the following loss function:
\begin{align}\label{eq:phiprime}
	\Psi&(\sigma,\mu) := \E \Big[ \textstyle \sum_{k=0}^{N} \theta_k \delta_k + \varsigma_k \Big],
\end{align}
for $\sigma$ that was obtained in the first step and for any $\mu \in \mathcal{M}$. Note that $\Psi(\sigma,\mu)$ is equivalent to $\Phi(\sigma,\mu)$. Associated with $\Psi(\sigma,\mu)$, we define the value function $V^c_k(\mathcal{I}^c_k)$ when $\sigma$ is given as
\begin{align}
	V^c_k(\mathcal{I}^c_k) :=& \min_{\mu \in \mathcal{M}}\E \Big[ \textstyle \sum_{t=k}^{N} \theta_{t-1} \delta_{t-1} + \varsigma_{t} \Big| \mathcal{I}^c_k \Big],\label{eq:Vc-def}
\end{align}
for $k \in \mathbb{N}_{[0,N]}$ with initial condition $V^c_{N+1}(\mathcal{I}^c_{N+1}) = 0$. By Lemmas~\ref{lemma:Kalmanfilter} and \ref{lemma:estimator-at-controller} in the Appendix, we observe that $\hat{e}_k$ and $\tilde{e}_k$ obey
\begin{align}
	\hat{e}_{k+1} &= A_k \hat{e}_k - \delta_k A_k \tilde{e}_k + w_k - (1-\delta_k) \imath_k, \label{eq:error-dyn}\\[2.25\jot]
	\tilde{e}_{k+1} &= (1-\delta_k) A_k \tilde{e}_k + K_{k+1} \nu_{k+1} - (1-\delta_k) \imath_k, \label{eq:mismatch-dyn}
\end{align}
for $k \in \mathbb{N}_{[0,N]}$ with initial conditions $\hat{e}_0 = x_0$ and $\tilde{e}_0 = K_0 \nu_0$, where $\imath_k = A_k \E[\hat{e}_k | \mathcal{I}^c_k,$ $\delta_k = 0]$. It is easy to deduce from (\ref{eq:error-dyn}) and (\ref{eq:mismatch-dyn}) that $\hat{e}_k$ and $\tilde{e}_k$ are independent of the control inputs under $\sigma$. Now, following a similar argument used in the proof of Theorem~1 in~\cite{voi}, we find that the value function $V^c_k(\mathcal{I}^c_k)$ should obey
\begin{align*}
	V^c_k(\mathcal{I}^c_k) &= \min_{u_k \in \mathbb{R}^m} \Big\{ \theta_{k-1} \E[ \delta_{k-1} | \mathcal{I}^c_k] + \tr(\Gamma_k Z_k)\\[1.5\jot]
	&\qquad \qquad + (u_k + L_k \hat{x}_k)^T (B_k^T S_{k+1} B_k + R_k)\\[2.25\jot]
	&\qquad \qquad \times  (u_k + L_k \hat{x}_k) + \E[V^c_{k+1}(\mathcal{I}^c_{k+1}) | \mathcal{I}^c_k] \Big\},
\end{align*}
for $k \in \mathbb{N}_{[0,N]}$, where $\delta_{k-1}$ and $Z_k = \Cov[ \hat{e}_k | \mathcal{I}^c_k]$ are independent of the control inputs. As a result, the minimizer is obtained by $u_k^\star = -L_k \hat{x}_k$. This establishes that $\Phi(\sigma,\xi) \leq \Phi(\sigma,\mu^o)$. 

In the third step, given the control policy $\xi$, we will prove that $\Phi(\omega,\xi) \leq \Phi(\sigma,\xi)$, where $\omega$ is a special form of $\sigma$ that is symmetric with respect to $\boldsymbol{\nu}_k$ at each time $k$. Let $\mathcal{N}$ be the set on which $\nu_k$ is defined, $\mathcal{B}(r)$ be a ball of radius $r$ centered at the origin and of proper dimension, and $\varpi_k \in \mathcal{N}$ be a variable obtained by the transformation $T_{k} \boldsymbol{\nu}_k$ for a given $T_{k}$. We recursively construct $\omega$ such that at each time $k$ the following conditions are satisfied:
\begin{equation}\label{eq:construction1}
\begin{aligned}
&\textstyle \int_{\mathcal{N}} \Prob_{\omega}(\delta_k = 0| \varpi_k,\boldsymbol{\delta}_{k-1} = 0) \ProbG_k(\varpi_k) d \varpi_k\\[2.25\jot]
&\ \ \qquad = \textstyle\int_{\mathcal{N}} \Prob_{\sigma}(\delta_k = 0| \varpi_k, \boldsymbol{\delta}_{k-1} = 0) \ProbQ_k(\varpi_k) d \varpi_k,
\end{aligned}
\end{equation}
and
\begin{equation}\label{eq:construction0}
\begin{aligned}
&\textstyle \int_{\mathcal{B}(r)} \Prob_{\omega}(\delta_{k}= 0| \varpi_k, \boldsymbol{\delta}_{k-1} = 0) \ProbG_k(\varpi_k) d \varpi_k \\[2.25\jot]
&\quad \textstyle \geq \int_{\mathcal{B}(r)} \big(\Prob_{\sigma}(\delta_{k}= 0| \varpi_k, \boldsymbol{\delta}_{k-1} = 0) \ProbQ_k(\varpi_k)\big)^* d \varpi_k,
\end{aligned}
\end{equation}
for all $r \geq 0$ with $ \Prob_{\omega}(\delta_k = 0| \varpi_k,\boldsymbol{\delta}_{k-1} = 0) \ProbG_k(\varpi_k)$ as a radially symmetric function of $\varpi_k$, where $\ProbG_k(\: . \:) := \Prob_{\omega}(\: . \: | \boldsymbol{\delta}_{k-1} = 0)$ and $\ProbQ_k(\: . \:) := \Prob_{\sigma}(\: . \: | \boldsymbol{\delta}_{k-1} = 0)$. Note that while the first condition states that $\Prob_{\omega}(\delta_{k}= 0|\varpi_k, \boldsymbol{\delta}_{k-1} = 0) \ProbG_k(\varpi_k)$ has the same volume under the curve as $(\Prob_{\sigma}(\delta_{k}= 0|\varpi_k, \boldsymbol{\delta}_{k-1} = 0) \ProbQ_k(\varpi_k))^*$, the second condition in fact states that the former is equally or more concentrated near the origin than the latter. This concentration near the origin, as we will see, leads to better estimation performance of the innovation, which is a Gaussian vector with zero mean.

Observe that
\begin{align*}
	\ProbG_{k+1}(\boldsymbol{\nu}_{k+1})& = \frac{\Prob(\nu_{k+1}) \Prob_{\omega}( \delta_{k} = 0 | \boldsymbol{\nu}_{k}, \boldsymbol{\delta}_{k-1} = 0) \ProbG_k(\boldsymbol{\nu}_k)}{\Prob_{\omega}(\delta_k = 0 | \boldsymbol{\delta}_{k-1} = 0)},\\[3.25\jot]
	\ProbQ_{k+1}(\boldsymbol{\nu}_{k+1}) &= \frac{\Prob(\nu_{k+1}) \Prob_{\sigma}( \delta_{k} = 0 | \boldsymbol{\nu}_{k}, \boldsymbol{\delta}_{k-1} = 0) \ProbQ_k(\boldsymbol{\nu}_k)}{\Prob_{\sigma}(\delta_k = 0 | \boldsymbol{\delta}_{k-1} = 0)},
\end{align*}
with initial conditions $\ProbG_0(\nu_0) = \ProbQ_0(\nu_0) = \Prob(\nu_0)$. Hence, given $T_{k}$, we can obtain $\ProbG_{k}(\varpi_{k})$ and $ \Prob_{\sigma}(\delta_k = 0| \varpi_k, \boldsymbol{\delta}_{k-1} = 0) \ProbQ_k(\varpi_k)$ based on $\ProbG_{k}(\boldsymbol{\nu}_{k})$ and $\ProbQ_{k+1}(\boldsymbol{\nu}_{k+1}) / \Prob(\nu_{k+1})$, respectively. Moreover, observe that
\begin{align*}
&\Prob_{\sigma}(\delta_k = 0 | \boldsymbol{\delta}_{k-1} = 0)\\[2.75\jot]
&=\textstyle \int_{\mathcal{N}} \Prob_{\sigma}(\delta_k = 0| \varpi_k, \boldsymbol{\delta}_{k-1} = 0) \Prob_{\sigma}(\varpi_k | \boldsymbol{\delta}_{k-1} = 0) d \varpi_k\\[2.75\jot]
&=\textstyle \int_{\mathcal{N}} \Prob_{\omega}(\delta_k = 0| \varpi_k, \boldsymbol{\delta}_{k-1} = 0) \Prob_{\omega}(\varpi_k | \boldsymbol{\delta}_{k-1} = 0) d \varpi_k\\[2.75\jot]
&=\Prob_{\omega}(\delta_k = 0 | \boldsymbol{\delta}_{k-1} = 0),
\end{align*}
where in the second equality we used (\ref{eq:construction1}). This relation will be useful in the following derivation.

To adopt the above construction, we need to introduce an equivalent loss function. It is possible to write
\begin{align*}
	\Psi(\sigma,\xi) &= \E \Big[ \textstyle \sum_{k=0}^{N} \theta_k \delta_k + \varsigma_k \Big]\\[1.5\jot]
	&= \E \Big[ \textstyle \sum_{k=0}^{N} \theta_k \delta_k + \hat{e}_k^T \Gamma_k \hat{e}_k \Big]\\[1.5\jot]
	& = \textstyle \sum_{k=0}^{N} \E \Big[ \theta_k \delta_k + \E [ \hat{e}_k^T \Gamma_k \hat{e}_k | \mathcal{I}^e_k] \Big]\\[1.5\jot]
	& = \textstyle \sum_{k=0}^{N} \E \Big[\theta_k \delta_k + \tilde{e}_k^T \Gamma_k \tilde{e}_k + \tr(\Gamma_k Y_k) \Big],
\end{align*}
for any $\sigma \in \mathcal{P}$ that is innovation-based and for $\xi$ that was obtained in the second step, where in the second equality we incorporated the control inputs $u_k = -L_k \hat{x}_k$, and in the third equality we used the tower property of conditional expectations. Note that $\Psi(\sigma,\xi)$ is equivalent to $\Phi(\sigma,\xi)$. Let us define the loss function $\Omega^M_{\sigma}(\tilde{e}_0)$~as
\begin{align*}
	\Omega^M_{\sigma}(\tilde{e}_0) := \textstyle \sum_{k=0}^{M} \E \Big[ \theta_k \delta_k + \tilde{e}_{k}^T \Gamma_{k} \tilde{e}_{k} \Big],
\end{align*}
for $M \in \mathbb{N}_{[0,N]}$ given $\tilde{e}_0$. Since $\tr(\Gamma_k Y_k)$ is independent of the decision policies, to prove the claim in the third step, it is enough to prove that $\Omega^M_{\omega}(\tilde{e}_0) \leq \Omega^M_{\sigma}(\tilde{e}_0)$ for any $M \in \{0,\dots,N\}$ and for any Gaussian vector $\tilde{e}_0$. Note that $\tilde{e}_0 = K_0 \nu_0$ under both $\sigma$ and $\omega$. Moreover, using the fact that $\Prob_{\sigma}(\delta_0 = 0) = \Prob_{\omega}(\delta_0 = 0)$, we~obtain
\begin{align*}
\E_{\sigma} \Big[ \delta_0 \Big] &=  1- \Prob_{\sigma}(\delta_0 = 0) \\[2.25\jot]
&= 1- \Prob_{\omega}(\delta_0 = 0) =\E_{\omega} \Big[ \delta_0 \Big] .
\end{align*}
Hence, the claim holds for the time horizon $0$. We assume that it also holds for all time horizons from $1$ to $M-1$. Observe that by the law of total probability, the following identities hold:
\begin{equation}\label{eq:identity2}
\begin{aligned}
	&\Prob_{\sigma}(\delta_0 = 1) + \Prob_{\sigma}(\boldsymbol{\delta}_{t} = 0)\\[1.75\jot]
	& \qquad \qquad  + \textstyle \sum_{s=1}^{t} \Prob_{\sigma}(\boldsymbol{\delta}_{s-1} = 0, \delta_s = 1) = 1,
\end{aligned}
\end{equation}
for any $t \in \mathbb{N}_{[0,N]}$. Applying the law of total expectation for the terms $\E[ \theta_k \delta_k]$ and $\E[\tilde{e}_{k}^T \Gamma_{k} \tilde{e}_{k}]$ in $\Omega^M_{\sigma}(\tilde{e}_0)$ on a partition provided by the identity (\ref{eq:identity2}) for $t = k-1$, and repeating this procedure for all $k \in \mathbb{N}_{[1,M]}$, we can obtain
\begin{align*}
	&\Omega^M_{\sigma}(\tilde{e}_0) =  \textstyle \sum_{k=0}^{M} \Big\{ \theta_k \Prob_{\sigma}(\boldsymbol{\delta}_{k-1} = 0) \E_{\sigma} \Big[ \delta_k \Big| \boldsymbol{\delta}_{k-1} = 0 \Big] \\[1.5\jot]
	&\qquad \qquad \qquad \qquad + \Prob_{\sigma}(\boldsymbol{\delta}_{k-1} = 0) \E_{\sigma} \Big[ \tilde{e}_{k}^T \Gamma_{k} \tilde{e}_{k} \Big| \boldsymbol{\delta}_{k-1} = 0 \Big]\\[2.25\jot]
	&\qquad \qquad \qquad \qquad + \Prob_{\sigma}(\boldsymbol{\delta}_{k-1} = 0, \delta_k = 1) \\[2.25\jot]
	&\qquad \qquad \qquad \qquad \times \E_{\sigma} \Big[ \Omega^{k+1,M}_{\sigma}(\tilde{e}_{k+1}) \Big| \boldsymbol{\delta}_{k-1} = 0, \delta_k = 1 \Big] \Big\},
\end{align*}
for $M \in \mathbb{N}_{[0,N]}$, where the cost-to-go $\Omega^{k,M}_{\sigma}(\tilde{e}_k)$ is defined as
\begin{align*}
	\Omega^{k,M}_{\sigma}(\tilde{e}_k) :=  \textstyle \sum_{t=k}^{M} \E \Big[ \theta_t \delta_t + \tilde{e}_{t}^T \Gamma_{t} \tilde{e}_{t} \Big],
\end{align*}
given $\tilde{e}_k$. Now, we will show that the probability coefficients, the transmission decision terms, the estimation mismatch terms, and the cost-to-go terms in $\Omega^M_{\sigma}(\tilde{e}_0)$ under $\sigma$ cannot be less than those when $\omega$ is used instead. First, note that since $\Prob_{\sigma}(\delta_k = 0 | \boldsymbol{\delta}_{k-1} = 0) = \Prob_{\omega}(\delta_k = 0 | \boldsymbol{\delta}_{k-1} = 0)$, we have $\Prob_{\sigma}(\boldsymbol{\delta}_{k-1} = 0) = \Prob_{\omega}(\boldsymbol{\delta}_{k-1} = 0) $ and $\Prob_{\sigma}(\boldsymbol{\delta}_{k-1} = 0, \delta_k = 1) = \Prob_{\omega}(\boldsymbol{\delta}_{k-1} = 0, \delta_k = 1)$. Hence, all the probability coefficients remain the same. Moreover, for the transmission decision terms, we~get
\begin{align*}
\E_{\sigma} \Big[ \delta_k \Big| \boldsymbol{\delta}_{k-1} = 0 \Big] &= 1- \Prob_{\sigma}(\delta_k = 0 | \boldsymbol{\delta}_{k-1} = 0)\\[2.25\jot]
& = 1- \Prob_{\omega}(\delta_k = 0 | \boldsymbol{\delta}_{k-1} = 0)\\[2.25\jot]
&=\E_{\omega} \Big[ \delta_k \Big| \boldsymbol{\delta}_{k-1} = 0 \Big].
\end{align*}
We continue the proof for the estimation mismatch terms by first showing that $\imath_k = 0$ for all $k \in \mathbb{N}_{[0,N]}$ under $\omega$. We assume that $\imath_t = 0$ for all $t \in \mathbb{N}_{[0,k-1]}$. It is possible to write
\begin{align*}
\E\Big[ \hat{e}_k \Big| \mathcal{I}_k^c, \delta_k \Big] &= \E \Big[ \E[\hat{e}_k | \mathcal{I}_k^e, \delta_k] \Big| \mathcal{I}_k^c, \delta_k \Big]\\[1\jot]
&= \E \Big[ \E[\hat{e}_k | \mathcal{I}_k^e ] \Big| \mathcal{I}_k^c, \delta_k \Big]\\[1\jot]
&=\E\Big[\tilde{e}_k \Big| \mathcal{I}_k^c, \delta_k \Big],
\end{align*}
where the first equality comes from the tower property of the conditional expectations and the second equality from the fact that $\delta_k$ is a function of $\mathcal{I}^e_k$. Hence, $\imath_k = A_k \E[\hat{e}_k | \mathcal{I}_k^c, \delta_k = 0 ] = A_k \E[\tilde{e}_k | \mathcal{I}_k^c, \delta_k = 0]$. Let $\tau_k$ denote the time elapsed since the last delivery when we are at time $k$. We have $\tilde{e}_{k-\tau_k} = K_{k-\tau_k} \nu_{k-\tau_k}$, and from (\ref{eq:mismatch-dyn}), we can express $\imath_k$ under $\omega$ as
\begin{align*}
	\imath_k &= A_k \E_{\omega} \Big[ \textstyle \sum_{t=0}^{\tau_k} D_{k-t} \nu_{k-t} \Big| \delta_{k-\tau_k} = 0,\dots, \delta_k = 0 \Big]\\[2\jot]
	&= A_k \textstyle \sum_{t=0}^{\tau_k} D_{k-t} \E_{\omega} \Big[ \nu_{k-t} \Big| \delta_{k-\tau_k} = 0,\dots, \delta_k = 0 \Big],
\end{align*}
where $D_{k-t}$ is a matrix depending on $A_{s}$ for $s \in \mathbb{N}_{[k-t, k-1]}$ and $K_{k-t}$. Since $\Prob_{\omega}(\boldsymbol{\nu}_k | \boldsymbol{\delta}_k = 0)$ has zero mean, we deduce that  $\Prob_{\omega}(\nu_{k-\tau_k}, \dots, \nu_k | \delta_{k-\tau_k} = 0,\dots, \delta_k = 0)$ has also zero mean. This implies that $\imath_k = 0$ for all $k \in \mathbb{N}_{[0,N]}$ under $\omega$. Given this observation, from (\ref{eq:mismatch-dyn}) when $\boldsymbol{\delta}_{k-1} = 0$, we find that $\tilde{e}_k = X_k \boldsymbol{\nu}_{k-1} + K_k \nu_k + c_k$ under $\sigma$, and that $\tilde{e}_k = X_k \boldsymbol{\nu}_{k-1} + K_k \nu_k$ under $\omega$, for a suitable matrix $X_k$ and a suitable vector $c_k$ both independent of $\boldsymbol{\nu}_k$. We can then write
\begin{align*}
&\E_{\sigma} \Big[ \tilde{e}_k^T \Gamma_k \tilde{e}_k \Big| \boldsymbol{\delta}_{k-1} = 0 \Big]\\[1.05\jot]
&\qquad = \E_{\sigma} \Big[ \big(X_k \boldsymbol{\nu}_{k-1} + K_k \nu_k + c_k \big)^T \Gamma_k\\[1.25\jot]
&\qquad \qquad \times \big(X_k \boldsymbol{\nu}_{k-1} + K_k \nu_k + c_k \big) \Big| \boldsymbol{\delta}_{k-1} = 0 \Big]\\[1.45\jot]
&\qquad = \E_{\sigma} \Big[ \boldsymbol{\nu}_{k-1}^T X_k^T \Gamma_k X_k \boldsymbol{\nu}_{k-1} + \nu_k^T K_k^T \Gamma_k K_k \nu_k \\[1.75\jot]
&\qquad \qquad + c_k^T \Gamma_k c_k + 2 \boldsymbol{\nu}_{k-1}^T X_k^T \Gamma_k c_k \Big| \boldsymbol{\delta}_{k-1} = 0 \Big],
\end{align*}
where in the second equality we used the fact that $\nu_k$ has zero mean and is independent of $\boldsymbol{\nu}_{k-1}$ and $\boldsymbol{\delta}_{k-1}$. Let us now use the decomposition $\Gamma_k = L_k^T U_k U_k^T L_k$, choose $T_{k-1} = U_k^T L_k X_k$, and define $f_{\sigma}(\varpi_{k-1}, \nu_k) := (\varpi_{k-1} + U_k^T L_k c_k)^T (\varpi_{k-1} + U_k^T L_k c_k) + \nu_k^T K_k^T \Gamma_k K_k \nu_k$, $f_{\omega}(\varpi_{k-1}, \nu_k) := \varpi_{k-1}^T \varpi_{k-1} + \nu_k^T K_k^T \Gamma_k K_k \nu_k$, $g_{\sigma}(\: . \:) := z -\min_z \{z,f_{\sigma}(\: . \:)\}$, and $g_{\omega}(\: . \:) := z -\min_z \{z,f_{\omega}(\: . \:)\}$. Clearly, for any fixed $z$, $g_{\sigma}(\varpi_{k-1}, \nu_k)$ and $g_{\omega}(\varpi_{k-1}, \nu_k)$ vanish at infinity. It follows that
\begin{align*}
	&\E_{\sigma} \Big[ \tilde{e}_k^T \Gamma_k \tilde{e}_k \Big| \boldsymbol{\delta}_{k-1} = 0 \Big] = \textstyle \int_{\mathcal{N}} \int_{\mathcal{N}} f_{\sigma}(\varpi_{k-1}, \nu_k) \\[2.75\jot]
	&\qquad \qquad \qquad \qquad \times  \Prob_{\sigma}(\varpi_{k-1} | \boldsymbol{\delta}_{k-1} = 0)  \Prob(\nu_{k}) d \varpi_{k-1} d\nu_k.
\end{align*}
In addition, we can write
\begin{align*}
	&\textstyle \int_{\mathcal{N}} g_{\sigma}(\varpi_{k-1}, \nu_k)\\[2.2\jot]
	&\quad \ \times \Prob_{\sigma}(\delta_{k-1} = 0| \varpi_{k-1}, \boldsymbol{\delta}_{k-2} = 0) \ProbQ_{k-1}(\varpi_{k-1}) d \varpi_{k-1}\\[2.5\jot]
	&\leq \textstyle  \int_{\mathcal{N}} g^*_{\sigma}(\varpi_{k-1}, \nu_k) \\[2.2\jot]
	&\quad \ \times \big( \Prob_{\sigma}(\delta_{k-1} = 0| \varpi_{k-1}, \boldsymbol{\delta}_{k-2} = 0) \ProbQ_{k-1}(\varpi_{k-1}) \big)^* d \varpi_{k-1}\\[2.2\jot]
	&= \textstyle \int_{\mathcal{N}} g_{\omega}(\varpi_{k-1}, \nu_k) \\[1.95\jot]
	&\quad \ \times \big( \Prob_{\sigma}(\delta_{k-1} = 0| \varpi_{k-1}, \boldsymbol{\delta}_{k-2} = 0) \ProbQ_{k-1}(\varpi_{k-1}) \big)^* d \varpi_{k-1}\\[2.45\jot]
	&\leq \textstyle \int_{\mathcal{N}} g_{\omega}(\varpi_{k-1}, \nu_k) \\[2.2\jot]
	&\quad \ \times \Prob_{\omega}(\delta_{k-1} = 0| \varpi_{k-1}, \boldsymbol{\delta}_{k-2} = 0) \ProbG_{k-1}(\varpi_{k-1}) d \varpi_{k-1},
\end{align*}
where in the first inequality we used the Hardy-Littlewood inequality (see Lemma~\ref{lemma:GHL} in the Appendix) with respect to $\varpi_{k-1}$, in the equality the fact that $g_{\sigma}^*(\varpi_{k-1}, \nu_k) = g_{\omega}(\varpi_{k-1}, \nu_k)$, and in the second inequality Lemma~\ref{lemma:major} in the Appendix and (\ref{eq:construction0}). This implies that
\begin{align*}
	&\textstyle \int_{\mathcal{N}} \textstyle \min_z \{z,f_{\sigma}(\varpi_{k-1}, \nu_k) \} \Prob_{\sigma}(\varpi_{k-1} | \boldsymbol{\delta}_{k-1} = 0) d \varpi_{k-1}\\[2.75\jot]
	&\geq \textstyle \int_{\mathcal{N}} \textstyle \min_z \{z,f_{\omega}(\varpi_{k-1}, \nu_k) \} \Prob_{\omega}(\varpi_{k-1} | \boldsymbol{\delta}_{k-1} = 0) d \varpi_{k-1},
\end{align*}
where we used the facts that
\begin{align*}
	&\Prob_{\sigma}(\varpi_{k-1} | \boldsymbol{\delta}_{k-1} = 0)\\[3.25\jot]
	&\qquad = \frac{\Prob_{\sigma}(\delta_{k-1} = 0 | \varpi_{k-1}, \boldsymbol{\delta}_{k-2} = 0) \ProbQ_{k-1}(\varpi_{k-1})}{\Prob_{\sigma}(\delta_{k-1} = 0|\boldsymbol{\delta}_{k-2} = 0)},
\end{align*}
and that $\Prob_{\sigma}(\delta_{k-1} = 0|\boldsymbol{\delta}_{k-2} = 0) = \Prob_{\omega}(\delta_{k-1} = 0|\boldsymbol{\delta}_{k-2} = 0)$. Now, taking $z$ to infinity, we conclude that
\begin{align*}
	&\textstyle \int_{\mathcal{N}} f_{\sigma}(\varpi_{k-1}, \nu_k) \Prob_{\sigma}(\varpi_{k-1} | \boldsymbol{\delta}_{k-1} = 0) d \varpi_{k-1}\\[2.25\jot]
	&\geq \textstyle \int_{\mathcal{N}} f_{\omega}(\varpi_{k-1}, \nu_k) \Prob_{\omega}(\varpi_{k-1} | \boldsymbol{\delta}_{k-1} = 0) d \varpi_{k-1}.
\end{align*}
Therefore,
\begin{align*}
\E_{\sigma} \Big[ \tilde{e}_k^T \Gamma_k \tilde{e}_k \Big| \boldsymbol{\delta}_{k-1} = 0 \Big] \geq \E_{\omega} \Big[ \tilde{e}_k^T \Gamma_k \tilde{e}_k \Big| \boldsymbol{\delta}_{k-1} = 0 \Big].	
\end{align*}
Finally, for the cost-to-go terms, we have 
\begin{align*}
	&\E_{\sigma} \Big[ \Omega^{k+1,M}_{\sigma}(\tilde{e}_{k+1}) \Big| \boldsymbol{\delta}_{k-1} = 0 , \delta_k = 1\Big]\\[2.25\jot]
	&=\textstyle \int_{\mathcal{N}^{k+2}}  \Omega^{k+1,M}_{\sigma}(\tilde{e}_{k+1}) \Prob_{\sigma}(\boldsymbol{\nu}_{k+1} | \boldsymbol{\delta}_{k-1} = 0, \delta_k = 1) d \boldsymbol{\nu}_{k+1}.
\end{align*}
Note that $\tilde{e}_{k+1} = K_{k+1} \nu_{k+1}$ under both $\sigma$ and $\omega$ when $\delta_k = 1$. Let $\bar{\Omega}^M_{\sigma}(\tilde{e}_0)$ denote a loss function that is structurally similar to $\Omega^M_{\sigma}(\tilde{e}_0)$ but with different parameters. Clearly, if $\Omega^M_{\sigma}(\tilde{e}_0) \geq \Omega^M_{\omega}(\tilde{e}_0)$, then $\bar{\Omega}^M_{\sigma}(\tilde{e}_0) \geq \bar{\Omega}^M_{\omega}(\tilde{e}_0)$. We can write
\begin{align*}
	&\textstyle \int_{\mathcal{N}^{k+2}} \Omega^{k+1,M}_{\sigma}(K_{k+1}\nu_{k+1})\\[1.85\jot]
	&\qquad \quad \times \Prob_{\sigma}(\boldsymbol{\nu}_{k+1} | \boldsymbol{\delta}_{k-1} = 0, \delta_k = 1) d \boldsymbol{\nu}_{k+1}\\[2.5\jot]
	&= \textstyle \int_{\mathcal{N}} \bar{\Omega}^{M-k-1}_{\sigma}(K_{k+1}\nu_{k+1}) \Prob(\nu_{k+1}) d \nu_{k+1}\\[2.75\jot]
	&\geq \textstyle \int_{\mathcal{N}} \bar{\Omega}^{M-k-1}_{\omega}(K_{k+1}\nu_{k+1}) \Prob(\nu_{k+1}) d \nu_{k+1}\\[2.5\jot]
	&= \textstyle  \int_{\mathcal{N}^{k+2}} \Omega^{k+1,M}_{\omega}(K_{k+1}\nu_{k+1}) \\[2.5\jot]
	&\qquad \quad \times \Prob_{\omega}(\boldsymbol{\nu}_{k+1} | \boldsymbol{\delta}_{k-1} = 0, \delta_k = 1) d \boldsymbol{\nu}_{k+1},
\end{align*}
where in the equalities we used the facts that $\Omega^{k+1,M}_{\sigma}(\tilde{e}) = \bar{\Omega}^{M-k-1}_{\sigma}(\tilde{e})$ for any Gaussian vector $\tilde{e}$ and a suitable selection of the parameters in $\bar{\Omega}^{M-k-1}_{\sigma}(\tilde{e})$, and that $\nu_{k+1}$ is independent of $\boldsymbol{\delta}_k$, and the Fubini's theorem; and in the inequality we used the hypothesis $\Omega^{M-k-1}_{\sigma}(\tilde{e}) \geq \Omega^{M-k-1}_{\omega}(\tilde{e})$ for any Gaussian vector $\tilde{e}$. Therefore,
\begin{align*}
&\E_{\sigma} \Big[ \Omega^{k+1,M}_{\sigma}(\tilde{e}_{k+1}) \Big| \boldsymbol{\delta}_{k-1} = 0 , \delta_k = 1\Big]\\[1.5\jot]
&\qquad \geq \E_{\omega} \Big[ \Omega^{k+1,M}_{\omega}(\tilde{e}_{k+1}) \Big| \boldsymbol{\delta}_{k-1} = 0 , \delta_k = 1\Big].
\end{align*}
This establishes that $\Omega^M_{\omega}(\tilde{e}_0) \leq \Omega^M_{\sigma}(\tilde{e}_0)$ and $\Phi(\omega,\xi) \leq \Phi(\sigma,\xi)$.

In the final step, we will conclude global optimality of the policy profile in the claim. Consider the following loss function:
\begin{align*}
	\Psi(\omega,\xi) &= \E \Big[ \textstyle \sum_{k=0}^{N} \theta_k \delta_k + \varsigma_k \Big],
\end{align*}
for any $\omega \in \mathcal{P}$ that is of the form specified in the third step and for $\xi$ that was obtained in the second step. Again note that $\Psi(\omega, \xi)$ is equivalent to $\Phi(\omega, \xi)$. Associated with $\Psi(\omega, \xi)$, we define the value function $V^e_k(\mathcal{I}^e_k)$ when $\xi$ is given as
\begin{align*}
	V^e_k(\mathcal{I}^e_k) := \min_{\omega \in \mathcal{P}}\E \Big[ \textstyle \sum_{t=k}^{N} \theta_t \delta_t + \varsigma_{t+1} \Big| \mathcal{I}^e_k \Big],
\end{align*}
for $k \in \mathbb{N}_{[0,N]}$ with initial condition $V^e_{N+1}(\mathcal{I}^e_{N+1}) = 0$ and with $\imath_t = 0$ for all $t \in \mathbb{N}_{[0,N]}$. Now, following a similar argument used in the proof of Theorem~1 in~\cite{voi}, we find that the value function $V^e_k(\mathcal{I}^e_k)$ should obey
\begin{align*}
	V^e_{k}(\mathcal{I}^e_k) &= \min_{\delta_k \in \{0,1\}} \Big\{\theta_k \delta_k + (1-\delta_k) \tilde{e}_k^T A_k^T \Gamma_{k+1} A_k \tilde{e}_k \\[1\jot]
	&\qquad \qquad \quad + \tr(A_k^T \Gamma_{k+1} A_k Y_k)\\[2.25\jot]
	&\qquad \qquad \quad +\tr(\Gamma_{k+1} W_k) + \E[V^e_{k+1}(\mathcal{I}^e_{k+1})|\mathcal{I}^e_k]  \Big\},
\end{align*}
for $k \in \mathbb{N}_{[0,N]}$. As a result, the minimizer is obtained by $\delta_k^\star = \mathds{1}_{\voi_k \geq0}$, where
\begin{align*}
	\voi_k &= \tilde{e}_k^T A_k^T \Gamma_{k+1} A_k \tilde{e}_k - \theta_k + \E[V^e_{k+1}(\mathcal{I}^e_{k+1})|\mathcal{I}^e_k, \delta_k = 0]\\[2.75\jot]
	&\quad - \E[V^e_{k+1}(\mathcal{I}^e_{k+1})|\mathcal{I}^e_k, \delta_k = 1].
\end{align*}
This certifies that $\Phi(\pi^\star,\mu^\star) \leq \Phi(\omega,\xi)$, and completes the proof.
\end{proof}

\section{Conclusion}\label{sec4}
In this article, we characterized a globally optimal solution in the rate-regulation tradeoff for multi-dimensional Gauss--Markov processes, and showed that such a solution consists of a symmetric threshold triggering policy based on the value of information and a certainty-equivalent control policy based on a non-Gaussian linear estimator. Besides, we provided a global optimality analysis for the value of information, and showed that it is globally optimal that the minimum mean-square-error state estimate at the event trigger or equivalently the estimation mismatch be transmitted to the controller only if the value of information becomes nonnegative. We suggest that future research should extend the framework developed in this study to more complex classes of systems.

\section*{Appendix}
In this section, we present a few lemmas that are used in our main analysis. The next two lemmas characterize the optimal estimators at the event trigger and the controller. For the proofs of these lemmas, see e.g., \cite{stengel1994} and \cite{voi}.
 
\begin{lemma}\label{lemma:Kalmanfilter}
The conditional mean $\E[{x}_k | \mathcal{I}^e_k]$ is the minimum mean-square-error estimator at the event trigger, and obeys
\begin{align}
\begin{split}\label{eq:kf-E}
	\check{x}_{k+1} &= A_k \check{x}_k + B_k u_k\\[1.75\jot]
	&\quad + K_{k+1} \big(y_{k+1} - C_{k+1} ( A_k \check{x}_k + B_k u_k)\big),
\end{split}\\[1.75\jot]
\begin{split}\label{eq:kf-Cov}
	Y_{k+1} &= \big( (A_k Y_k A_k^T + W_k)^{-1} + C_{k+1}^T V_{k+1}^{-1} C_{k+1} \big)^{-1},	
\end{split}
\end{align}
for $k \in \mathbb{N}_{[0,N]}$ with initial conditions $\check{x}_0 = m_0 + Y_{0} C_{0}^T V_{0}^{-1}(y_0 - C_0 m_0)$ and $Y_0 = (M_0^{-1} + C_{0}^T V_{0}^{-1} C_{0})^{-1}$, where $\check{x}_k = \E[{x}_k | \mathcal{I}^e_k]$, $Y_k = \Cov[x_k | \mathcal{I}^e_k]$, and $K_{k} = Y_{k} C_{k}^T V_{k}^{-1}$.
\end{lemma}

\begin{lemma}\label{lemma:estimator-at-controller}
The conditional mean $\E[{x}_k | \mathcal{I}^c_k]$ is the minimum mean-square-error estimator at the controller, and obeys
\begin{equation}\label{c1:eq:imperf-estimate-dynX}
\begin{aligned}
	\hat{x}_{k+1} &= A_k \hat{x}_k + B_k u_k + \delta_k A_k \tilde{e}_k + (1-\delta_k) \imath_k,
\end{aligned}
\end{equation}
for $k \in \mathbb{N}_{[0,N]}$ with initial condition $\hat{x}_0 = m_0$, where $\hat{x}_k = \E[{x}_k | \mathcal{I}^c_k]$ and $\imath_k = A_k \E[\hat{e}_k | \mathcal{I}^c_k,\delta_k=0]$. In addition, the conditional covariance $\Cov[{x}_k | \mathcal{I}^c_k]$ obeys
\begin{equation}\label{c1:eq:imperf-cov-dynX}
\begin{aligned}
	Z_{k+1} &= A_k Z_k A_k^T + W_k\\[1.75\jot]
	&\qquad - \delta_k A_k (Z_k - Y_k) A_k^T - (1-\delta_k) \Xi_k,	
\end{aligned}
\end{equation}
for $k \in \mathbb{N}_{[0,N]}$ with initial condition $Z_0 = M_0$, where $Z_k = \Cov[x_k | \mathcal{I}^c_k]$ and $\Xi_k = A_k (Z_k - \Cov[\hat{e}_k | \mathcal{I}_k^c, \delta_k =0]) A_k^T$.
\end{lemma}

Moreover, the next two lemmas are pertaining to symmetric decreasing rearrangements of non-negative functions. For the proofs of these lemmas, see e.g., \cite{brock2000} and \cite{alvino1991}.

\begin{lemma}[Hardy-Littlewood inequality]\label{lemma:GHL}
Let $f$ and $g$ be non-negative functions defined on $\mathbb{R}^n$ that vanish at infinity. Then,
\begin{align}
	\textstyle \int_{\mathbb{R}^n} f(x) g(x) dx \leq \int_{\mathbb{R}^n} f^*(x) g^*(x) dx.
\end{align}
\end{lemma}
\vspace{2mm}

\begin{lemma}\label{lemma:major}
Let $\mathcal{B}(r) \subseteq \mathbb{R}^n$ be a ball of radius $r$ centered at the origin, and $f$ and $g$ be non-negative functions defined on $\mathbb{R}^n$ that vanish at infinity and obey
\begin{align}
	\textstyle \int_{\mathcal{B}(r)} f^*(x) dx \leq \int_{\mathcal{B}(r)} g^*(x) dx,
\end{align}
for all $r \geq 0$. Then,
\begin{align}
	\textstyle \int_{\mathcal{B}(r)} h(x) f^*(x) dx \leq \int_{\mathcal{B}(r)} h(x) g^*(x) dx,
\end{align}
for any symmetric non-increasing function $h$.
\end{lemma}

\bibliography{../../mybib}
\bibliographystyle{ieeetr}

\end{document}